\documentclass[11pt]{article}

\usepackage{amssymb}
\usepackage{amsmath}
\usepackage{color}
\usepackage{theorem}
\usepackage[all]{xy}

\newcommand{\XYMATRIX}{\xymatrix@M=6pt}
\newcommand{\XYMATRIXNARROW}{\xymatrix@C=2pt@R=6pt}
\newcommand{\aremb}{\ar@{^{(}->}}
\newcommand{\arembfrom}{\ar@{<-^{)}}}

\numberwithin{equation}{section}

\theorembodyfont{\slshape}

  \newtheorem{THM}{Theorem}[section]

  \newtheorem{COR}[THM]{Corollary}

\theorembodyfont{\rmfamily}

  \newtheorem{DEF}[THM]{Definition}
  \newtheorem{EX}{Example}[section]

\newif\ifQEDsign
\newcommand{\QED}{\global\QEDsigntrue\hfill$\square$}

\newenvironment{proof}%
    {\par\noindent\textit{Proof.}\global\QEDsignfalse}%
    {\ifQEDsign\else\QED\fi\par\bigskip\par}

\def\labelenumi{(\roman{enumi})}

\renewcommand{\le}{\leqslant}
\renewcommand{\ge}{\geqslant}

\newcommand{\0}{\varnothing}

\renewcommand{\sec}{\cap}
\renewcommand{\phi}{\varphi}
\renewcommand{\epsilon}{\varepsilon}

\newcommand{\CC}{\mathbf{C}}
\newcommand{\DD}{\mathbf{D}}

\newcommand{\KK}{\mathbf{K}}

\newcommand{\NN}{\mathbb{N}}

\newcommand{\union}{\cup}
\newcommand{\restr}[2]{\hbox{$#1$}\hbox{$\upharpoonright$}_{#2}}
\newcommand{\reduct}[2]{\hbox{$#1$}\hbox{$|$}_{#2}}
\newcommand{\Boxed}[1]{\mbox{$#1$}}

\newcommand{\id}{\mathrm{id}}

\newcommand{\Ob}{\mathrm{Ob}}

\newcommand{\arity}{\mathrm{ar}}

\newcommand{\calA}{\mathcal{A}}
\newcommand{\calB}{\mathcal{B}}
\newcommand{\calC}{\mathcal{C}}
\newcommand{\calD}{\mathcal{D}}

\newcommand{\calS}{\mathcal{S}}
\newcommand{\calT}{\mathcal{T}}

\newcommand{\calX}{\mathcal{X}}

\newcommand{\CHemb}{\mathbf{Ch}}
\newcommand{\EPosemb}{\mathbf{EPos}}

\title{A Ramsey Theorem for Multiposets}
\author{%
  Nemanja Dragani\'c\\
  Department of Mathematics\\
  ETH Z\"urich, Z\"urich, Switzerland\\
  email: dnemanja@ethz.ch
\and  
  Dragan Ma\v sulovi\'c\\
  University of Novi Sad, Faculty of Sciences\\
  Department of Mathematics and Informatics\\
  Trg Dositeja Obradovi\'ca 3, 21000 Novi Sad, Serbia\\
  e-mail: dragan.masulovic@dmi.uns.ac.rs}
  
\begin{document}
\maketitle

\begin{abstract}
  In the parlance of relational structures, the Finite Ramsey Theorem states that the class of all finite
  chains has the Ramsey property. A classical result from the 1980's claims that the class of all finite
  posets with a linear extension has the Ramsey property. In 2010 Soki\'c proved that the
  class of all finite structures consisting of several linear orders has the Ramsey property.
  This was followed by a 2017 result of Solecki and Zhao that the class of all finite
  posets with several linear extensions has the Ramsey property.

  Using the categorical reinterpretation of the Ramsey property in this paper we prove
  a common generalization of all these results. We consider multiposets to be structures consisting of
  several partial orders and several linear orders. We allow partial orders to extend each other
  in an arbitrary but fixed way, and require that every partial order is extended by at least one of the linear orders.
  We then show that the class of all finite multiposets conforming to a fixed template has the Ramsey property.

  \bigskip


  \noindent \textbf{AMS Subj.\ Classification (2010):} 05C55, 18A99
\end{abstract}

\section{Introduction}

Generalizing the classical results of F.~P.~Ramsey from the late 1920's, the structural Ramsey theory originated at
the beginning of 1970's in a series of papers (see \cite{N1995} for references).
We say that a class $\KK$ of finite structures has the \emph{Ramsey property} if the following holds:
for any number $k \ge 2$ of colors and all $\calA, \calB \in \KK$ such that $\calA$ embeds into $\calB$
there is a $\calC \in \KK$
such that no matter how we color the copies of $\calA$ in $\calC$ with $k$ colors, there is a \emph{monochromatic} copy
$\calB'$ of $\calB$ in $\calC$ (that is, all the copies of $\calA$ that fall within $\calB'$ are colored by the same color).
In this parlance, the Finite Ramsey Theorem~\cite{Ramsey} takes the following form:

\begin{quote}
  \textsl{%
    (Finite Ramsey Theorem) The class of all finite cha\-ins has the Ramsey property.
  }
\end{quote}

As it turns out many classes of finite linearly ordered structures have the Ramsey property.
For example, the class of all finite linearly ordered graphs has the Ramsey property~\cite{AH,Nesetril-Rodl}.
Interestingly, this is not the case with the class of finite partial orders accompanied with
arbitrary linear orders --- this class is not Ramsey~\cite{fouche,sokic-phd}. However,
Paoli, Trotter and Walker show in~\cite{PTW} that the class of all finite posets with a linear extension has
the Ramsey property.\footnote{%
  Let us briefly note that the proof of the ordering property presented in
  \cite[Theorem~16, p.~362]{PTW} relies on Theorem~2 (p.~354), the statement of which is not true.
  Fortunately,
  as the referess of this paper have pointed out to the second author, the finitary version
  of Theorem~2 (the well known Finite Product Ramsey Theorem) suffices for the proof.
}
The same result was reproved recently using different strategies by
Soki\'c~\cite{sokic-phd,sokic-order}, Solecki and Zhao~\cite{solecki-zhao},
Ne\v set\v ril and R\"odl~\cite{Nesetril-Rodl-2016}, and the second author~\cite{masul-preadj}.
Soki\'c in~\cite{sokic-phd,sokic-order} derived the result as a consequence of a result of Fouch\'e~\cite{fouche},
while the proof of Solecki and Zhao~\cite{solecki-zhao} relies on Solecki's
abstract approach to finite Ramsey theory~\cite{solecki-abstract}.
The proof given in~\cite{Nesetril-Rodl-2016} starts from
the fact that the class of all finite acyclic digraphs endowed with a linear extension is a Ramsey class
(and this, in turn, follows from the Ne\v set\v ril-R\"odl Theorem) and then uses the partite construction
to ``improve acyclic digraphs to posets''. This is the first proof where the partite construction was used to establish
the Ramsey property for this class of structures.
The proof given in~\cite{masul-preadj}, on the other hand,
establishes a particular relationship (called a pre-adjunction) between the category of all finite posets with a linear extension
and the Graham-Rothschild category (which is nothing but a categorical rendering of the setup of the
Graham-Rothschild Theorem) and then uses the pre-adjunction between the categories to ``transport the Ramsey property''
from the the Graham-Rothschild category to the other one. Interestingly, the proofs presented in
\cite{solecki-zhao} and \cite{masul-preadj} do not require the ordering property in order to
establish the Ramsey property for the class.

The fact that the class of all finite posets with a linear extension has the Ramsey property
was generalized in several directions. In his paper ~\cite{fouche} Fouch\'e explicitly calculated
Ramsey degrees of finite posets: as it turns out, the Ramsey degree of a finite
poset is the number of finite linear extensions of the poset ordering. Generalizing the result in
other direction, Soki\'c proved in his PhD thesis~\cite[Theorem~78, p.~96]{sokic-phd} that the
class of all finite structures consisting of several linear orders has the Ramsey property.
(The case for $n = 2$ was independently proved by B\"ottcher and Foniok in~\cite{bottcher-foniok}.)
Soki\'c also proved that the class of all finite posets with a linear extension and an independent linear order has the Ramsey
property~\cite[Theorem~80, p.~98]{sokic-phd}.
This example is interesting and motivating, because it seems that as soon as we have at least one linear extension
of the base partial order we can pretty much do whatever we want.

Another point of view on the same phenomenon is taken by Soki\'c in~\cite{sokic-phd} and
Bodirsky in~\cite{bodirsky1,bodirsky2}. They proved that
we can always ``put two Ramsey classes together'' to get a new one
in the following sense: free interposition of two Ramsey classes with strong amalgamation
is again a Ramsey class.

The next important step in understanding the Ramsey property for classes of finite partial orders endowed with
additional linear orders was a result of Solecki and Zhao that the class of all finite
posets with several linear extensions has the Ramsey property~\cite{solecki-zhao}.
An alternative, shorter, proof of the Solecki-Zhao result was given by Arman and R\"odl in~\cite{arman-rodl}.
The proof of Arman and R\"odl starts from a large product of classes of finite linear orders with
a linear extension and then reduces it to the class of finite posets with several linear extensions.
In this paper we take a similar approach, but in order to prove our more general result we have to refine it.
As in the case of the proof of Arman and R\"odl our main ``building tool'' is the structural version of the Product
Ramsey Theorem from~\cite{sokic-boron}, but as the ``refinement tool'' we use a theorem about transferring the Ramsey property from a
category onto its subcategory closed in a particular way.

A problem closely related to identifying Ramsey classes is the classification of amalgamation classes of finite
structures (a connection between the two notions was established by Ne\v set\v ril in~\cite{nesetril}).
Amalgamation classes of permutations (understood as finite structures with two independent linear orders) were classified by
Cameron in~\cite{cameron-perm}. In that paper Cameron also posed the problem of generalizing his result
to classes of finite structures with three or more independent linear orders, which turned out to be
quite a challenge. A very recent result of Braunfeld and Simon \cite{braunfeld-simon} gives a catalog of amalgamation
classes of such structures which, interestingly, contains examples of classes of structures that do not fall into
the framework of this paper. The relationship of the classes from their catalog to the Ramsey property
(in the context of Ramsey expansions) was further discussed by Braunfeld in~\cite{braunfeld}.

Going back to the examples that motivate the main result of this paper let us collect
some of the above considerations into a single statement.

\begin{THM}\label{mp.thm.ALL}
  \begin{enumerate}\renewcommand{\labelenumi}{$(\alph{enumi})$}
  \item
    (Ramsey~\cite{Ramsey})
    The class of all finite cha\-ins $(A, \Boxed\le)$ has the Ramsey property.
  \item
    (Paoli, Trotter and Walker~\cite{PTW})
    The class of all finite posets with a linear extension (that is, structures $(A, \Boxed{\le_1}, \Boxed{\le_2})$ where
    $\le_1$ is a partial order on $A$ and $\le_2$ a linear order on $A$ such that $(\Boxed{\le_1}) \subseteq (\Boxed{\le_2})$)
    has the Ramsey property.
  \item
    (Soki\'c~\cite{sokic-phd})
    The class of all finite posets with a linear extension and another linear order (that is, structures $(A, \Boxed{\le_1}, \Boxed{\le_2}, \Boxed{\le_3})$
    where $\le_1$ is a partial order, $\le_2$ a linear extension of~$\le_1$ and $\le_3$ is an arbitrary linear order)
    has the Ramsey property.
  \item
    (Soki\'c~\cite{sokic-phd,sokic-order}, Bodirsky~\cite{bodirsky1,bodirsky2}, B\"ottcher and Foniok~\cite{bottcher-foniok})
    For every $n \ge 2$, the class of
    all finite structures of the form $(A, \Boxed{\le_1}, \ldots, \Boxed{\le_n})$ where each $\le_i$
    is a linear order on $A$ has the Ramsey property.
  \item
    (Solecki and Zhao~\cite{solecki-zhao}, Arman and R\"odl~\cite{arman-rodl})
    For every $n \ge 1$, the class of
    all finite structures of the form $(A, \Boxed{\le_0}, \Boxed{\le_1}, \ldots, \Boxed{\le_n})$ where $\le_0$
    is a partial order on $A$ and each $\le_i$, $i \in \{1, 2, \ldots, n\}$, is a linear order on $A$ extending $\le_0$
    has the Ramsey property.
  \end{enumerate}
\end{THM}

\begin{figure}
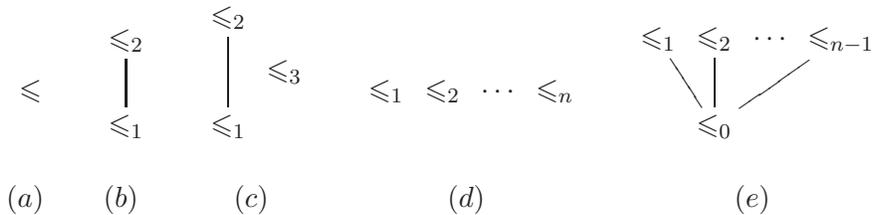

  \centering
  \begin{tabular}[b]{c}
    $$
      \XYMATRIXNARROW{
        \\
        \Boxed{\le} \\
        \\
      }
    $$
    \\
    \ \\
    $(a)$
  \end{tabular}
  \begin{tabular}[b]{c}
    $$
      \XYMATRIXNARROW{
        \Boxed{\le_2} \ar@{-}[dd] \\
        \\
        \Boxed{\le_1}\\
      }
    $$
    \\
    \ \\
    $(b)$
  \end{tabular}
  \begin{tabular}[b]{c}
    $$
      \XYMATRIXNARROW{
        \Boxed{\le_2} \ar@{-}[dd] &               \\
                                  & \Boxed{\le_3} \\
        \Boxed{\le_1}             &               \\
      }
    $$
    \\
    \ \\
    $(c)$
  \end{tabular}
  \begin{tabular}[b]{c}
    $$
      \XYMATRIXNARROW{
                       &               &        &                \\
         \Boxed{\le_1} & \Boxed{\le_2} & \cdots &  \Boxed{\le_n} \\
                       &               &        &                \\
      }
    $$
    \\
    \ \\
    $(d)$
  \end{tabular}
  \begin{tabular}[b]{c}
    $$
      \XYMATRIXNARROW{
          \Boxed{\le_1} \ar@{-}[ddr] & \Boxed{\le_2} \ar@{-}[dd] & \cdots & \Boxed{\le_{n-1}} \ar@{-}[ddll] \\
          \\
                                     & \Boxed{\le_0}\\
      }
    $$
    \\
    \ \\
    $(e)$
  \end{tabular}
\caption{Relationships between ordering relations in five situations listed in Theorem~\ref{mp.thm.ALL}}
\label{mp.fig.1}
\end{figure}

For reasons that will become apparent soon, let us depict these five situations as in Fig.~\ref{mp.fig.1}.
For example, by Fig.~\ref{mp.fig.1}~$(c)$ we indicate that the structures we are interested in have three ordering relations,
the first one is always contained in the second, while the third ordering relation is independent of the first two.

In each of these cases we have a class of structures with several ordering relations, the relations are required to form a
fixed partially ordered set under set inclusion, and the maximal elements in this poset of relations are required to
be linear orders. A straightforward generalization now leads to the following concept.

Let $\calT = (\{1, 2, \ldots, t\}, \Boxed\preccurlyeq)$, $t \ge 1$, be a poset which we refer to as the \emph{template}.
A \emph{$\calT$-multiposet} is a structure $(A, \Boxed{\le_1}, \ldots, \Boxed{\le_t})$ where
\begin{itemize}
\item
  $\Boxed{\le_1}$, \ldots, $\Boxed{\le_t}$ are partial orders on $A$,
\item
  if $i$ is a maximal element of $\calT$ then $\le_i$ is a linear order on $A$, and
\item
  if $i \preccurlyeq j$ in $\calT$ then $(\Boxed{\le_i}) \subseteq (\Boxed{\le_j})$.
\end{itemize}
Let $\KK(\calT)$ be the class of all finite $\calT$-multiposets.
The purpose of this paper is to show the following result which
clearly generalizes each of the results listed in Theorem~\ref{mp.thm.ALL}:

\begin{THM}\label{mp.thm.MAIN}
  For every template $\calT$ the class $\KK(\calT)$ has the Ramsey  property.
\end{THM}

It was Leeb who pointed out in 1970 \cite{leeb-cat} that the use of category theory can be quite helpful
both in the formulation and in the proofs of results pertaining to structural Ramsey theory.
We pursued this line of thought in several papers \cite{masulovic-ramsey,masul-preadj,masul-drp-perm}
and demonstrated that reinterpreting the Ramsey property in the context of category theory and using the
machinery of category theory can lead to essentially new proving strategies. The proof of
Theorem~\ref{mp.thm.MAIN} will represent another demonstration of these new
strategies.

In Section~\ref{mp.sec.prelim} we
give a brief overview of standard notions referring to finite structures and category theory,
and conclude with the reinterpretation of the Ramsey property in the language of category theory.
In Section~\ref{mp.sec.transfer} we consider two ways of transferring the Ramsey property from a category
to another category. We first recall a result of M.~Soki\'c from~\cite{sokic-boron}
which enables us to combine Ramsey classes of finite structures over disjoint relational signatures in a particular way,
and then recall a result from~\cite{masul-drp-perm} which enables us to transfer the Ramsey property from a category
to its (not necessarily full) subcategory. Using these two transfer principles, starting from
Theorem~\ref{mp.thm.ALL}~$(a)$ and~$(b)$, in Section~\ref{mp.sec.proof} we prove Theorem~\ref{mp.thm.MAIN}.

\section{Preliminaries}
\label{mp.sec.prelim}

In this section we give a brief overview of standard notions referring to first order structures
and conclude with a reinterpretation of the Ramsey property in the language of category theory.

\subsection{Structures}

Let $\Theta$ be a set of relational symbols.
A \emph{$\Theta$-structure} $\calA = (A, \Theta^\calA)$ is a set $A$ together with a set $\Theta^\calA$ of
relations on $A$ which are interpretations of the corresponding symbols in $\Theta$.
A structure $\calA = (A, \Theta^\calA)$ is \emph{finite} if $A$ is a finite set.
For $\Theta$-structures $\calA$ and $\calB$,
a \emph{embedding} $f: \calA \hookrightarrow \calB$ is an injection $f: A \rightarrow B$ such that
$(a_1, \ldots, a_r) \in \rho^\calA \Leftrightarrow (f(a_1), \ldots, f(a_r)) \in \rho^\calB$,
for every relational symbol $\rho \in \Theta$ and all $a_1, \ldots, a_r \in A$ where $r = \arity(\rho)$.

A structure $\calA$ is a \emph{substructure} of a structure
$\calB$ ($\calA \le \calB$) if the identity map $a \mapsto a$ is an embedding of $\calA$ into $\calB$.
Let $\calA$ be a structure and $\0 \ne B \subseteq A$. Then $\restr \calA B = (B, \restr{\Theta^\calA}{B})$ denotes
the \emph{substructure of $\calA$ induced by~$B$}, where $\restr{\Theta^\calA}{B}$ denotes the restriction of each
relation in $\Theta^\calA$ to~$B$.
If $\calA$ is a $\Theta$-structure and $\Sigma \subseteq \Theta$ then by $\reduct \calA \Sigma$ we denote
the \emph{$\Sigma$-reduct} of~$\calA$: $\reduct \calA \Sigma = (A, \{\theta^\calA : \theta \in \Sigma\})$.

Let $\Theta$ be a set of relational symbols.
Let $\CC$ be a class of $\Theta$-structures and $\KK$ a subclass of $\CC$.
Then:
\begin{itemize}
\item
  $\KK$ has the \emph{hereditary property (HP) with respect to $\CC$} if
  for all $\calC \in \KK$ and $\calB \in \CC$ such that $\calB \hookrightarrow \calC$ we have that $\calB \in \KK$;
\item
  $\KK$ has the \emph{joint embedding property (JEP)} if for all $\calA, \calB \in \KK$ there is a $\calC \in \KK$
  such that $\calA \hookrightarrow \calC$ and $\calB \hookrightarrow \calC$;
\item
  $\KK$ has the \emph{strong amalgamation property (SAP)} if for all $\calA, \calB, \calC \in \KK$ and embeddings
  $f_1 : \calA \hookrightarrow \calB$ and $f_2 : \calA \hookrightarrow \calC$ there is a $\calD \in \KK$ and embeddings
  $g_1 : \calB \hookrightarrow \calD$ and $g_2 : \calC \hookrightarrow \calD$ such that $g_1 \circ f_1 = g_2 \circ f_2$ and
  $g_1(B) \sec g_2(C) = g_1 \circ f_1 (A) = g_2 \circ f_2 (A)$.
\end{itemize}

\begin{EX}\label{mp.ex.class-HP-JEP-SAP}
  Let $\CHemb$ denote the class of all finite chains, and let $\EPosemb$ denote the class of all finite posets with
  a linear extension. Both $\CHemb$ and $\EPosemb$ have each of the properties listed above.
\end{EX}

\subsection{Categories, functors and the Ramsey property}

In order to specify a \emph{category} $\CC$ one has to specify
a class of objects $\Ob(\CC)$, a set of morphisms $\hom_\CC(\calA, \calB)$ for all $\calA, \calB \in \Ob(\CC)$,
the identity morphism $\id_\calA$ for all $\calA \in \Ob(\CC)$, and
the composition of mor\-phi\-sms~$\cdot$~so that
$\id_\calB \cdot f = f \cdot \id_\calA$ for all $f \in \hom_\CC(\calA, \calB)$, and
$(f \cdot g) \cdot h = f \cdot (g \cdot h)$.
A morphism $f \in \hom_\CC(B, C)$ is \emph{monic} or \emph{left cancellable} if
$f \cdot g = f \cdot h$ implies $g = h$ for all $g, h \in \hom_\CC(A, B)$ where $A \in \Ob(\CC)$ is arbitrary.

\begin{EX}
  Every class of structures of the same relational type forms a category where morphisms are embeddings.
  Given a class $\KK$ of structures (of the same relational type), whenever we refer to $\KK$ as a category,
  we have in mind the category whose objects are structures from $\KK$ and whose morphisms are embeddings.
  
  In particular, $\CHemb$ and $\EPosemb$ are categories where objects are the corresponding structures
  and morphisms are embeddings.
\end{EX}

A category $\DD$ is a \emph{subcategory} of a category $\CC$ if $\Ob(\DD) \subseteq \Ob(\CC)$ and
$\hom_\DD(\calA, \calB) \subseteq \hom_\CC(\calA, \calB)$ for all $\calA, \calB \in \Ob(\DD)$.
A category $\DD$ is a \emph{full subcategory} of a category $\CC$ if $\Ob(\DD) \subseteq \Ob(\CC)$ and
$\hom_\DD(\calA, \calB) = \hom_\CC(\calA, \calB)$ for all $\calA, \calB \in \Ob(\DD)$.

A \emph{functor} $F : \CC \to \DD$ from a category $\CC$ to a category $\DD$ maps $\Ob(\CC)$ to
$\Ob(\DD)$ and maps morphisms of $\CC$ to morphisms of $\DD$ so that
$F(f) \in \hom_\DD(F(\calA), F(\calB))$ whenever $f \in \hom_\CC(\calA, \calB)$, $F(f \cdot g) = F(f) \cdot F(g)$ whenever
$f \cdot g$ is defined, and $F(\id_\calA) = \id_{F(\calA)}$.

Categories $\CC$ and $\DD$ are \emph{isomorphic} if there exist functors $F : \CC \to \DD$ and $G : \DD \to \CC$ which are
inverses of one another both on objects and on morphisms.

\begin{figure}[h]
  $$
  \xymatrix{
    & & & & & \exists C &
  \\
    \bullet & \bullet & \bullet
    & & B_1 \ar@{.>}[ur] & B_2 \ar@{.>}[u] & B_3 \ar@{.>}[ul]
  \\
    \bullet \ar[u] \ar[ur] & \bullet \ar[ur] \ar[ul] & \bullet \ar[ul] \ar[u]
    & & A_1 \ar[u]^{f_1} \ar[ur]_(0.3){f_2} & A_2 \ar[ur]^(0.3){f_4} \ar[ul]_(0.3){f_3} & A_3 \ar[ul]^(0.3){f_5} \ar[u]_{f_6}
  \\
    & \Delta \ar[rrrr]^F  & & & & \CC  
  }
  $$
  \caption{A consistent diagram in $\CC$ (of shape $\Delta$)}
  \label{catsap.fig.3}
\end{figure}
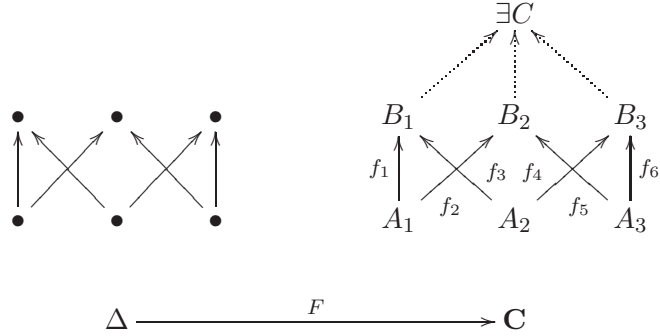

A \emph{diagram} in a category $\CC$ is a functor $F : \Delta \to \CC$ where the category $\Delta$ is referred to as the
\emph{shape of the diagram}. A diagram $F : \Delta \to \CC$ is \emph{consistent in $\CC$} if there exists a $C \in \Ob(\CC)$
and a family of morphisms $(e_\delta : F(\delta) \to C)_{\delta \in \Ob(\Delta)}$ such that for every
morphism $g : \delta \to \gamma$ in $\Delta$ we have $e_\gamma \cdot F(g) = e_\delta$:
$$
  \xymatrix{
     & C & \\
    F(\delta) \ar[ur]^{e_\delta} \ar[rr]_{F(g)} & & F(\gamma) \ar[ul]_{e_\gamma}
  }
$$
(see Fig.~\ref{catsap.fig.3}). We say that $C$ together with the family of morphisms
$(e_\delta)_{\delta \in \Ob(\Delta)}$ forms a \emph{compatible cone in $\CC$ over the diagram~$F$}.

Let $\CC$ be a category and $\calS$ a set. We say that
$
  \calS = \calX_1 \union \ldots \union \calX_k
$
is a \emph{$k$-coloring} of $\calS$ if $\calX_i \sec \calX_j = \0$ whenever $i \ne j$.
For an integer $k \ge 2$ and $\calA, \calB, \calC \in \Ob(\CC)$ we write
$
  \calC \longrightarrow (\calB)^{\calA}_k
$
to denote that for every $k$-coloring
$
  \hom_\CC(\calA, \calC) = \calX_1 \union \ldots \union \calX_k
$
there is an $i \in \{1, \ldots, k\}$ and a morphism $w \in \hom_\CC(\calB, \calC)$ such that
$w \cdot \hom_\CC(\calA, \calB) \subseteq \calX_i$.

\begin{DEF}
  A category $\CC$ has the \emph{Ramsey property} if
  for every integer $k \ge 2$ and all $\calA, \calB \in \Ob(\CC)$
  such that $\hom_\CC(\calA, \calB) \ne \0$ there is a
  $\calC \in \Ob(\CC)$ such that $\calC \longrightarrow (\calB)^{\calA}_k$.
\end{DEF}

Clearly, if $\CC$ and $\DD$ are isomorphic categories, then one of them has the Ramsey property if and only if
the other one does.

\begin{EX}\label{mp.ex.class-RP}
    As we have seen, both $\CHemb$ and $\EPosemb$ have the Ramsey property
    (Theorem~\ref{mp.thm.ALL}~$(a)$ and~$(b)$).
\end{EX}

\section{Transferring the Ramsey property between categories}
\label{mp.sec.transfer}

In this section we give a brief overview of two strategies of transferring the Ramsey property from a category
to another category. We first recall a result of M.~Soki\'c from~\cite{sokic-boron}
which enables us to combine Ramsey classes of structures over disjoint signatures in a particular way.

Let $\Sigma_1$, $\Sigma_2$, \ldots, $\Sigma_n$, $n \ge 2$, be pairwise disjoint sets of relational symbols,
and for each $i \in \{1, 2, \ldots, n\}$ let $\KK_i$ be a class of $\Sigma_i$-structures. Let $\Theta = \bigcup_{i=1}^n \Sigma_i$.
Then with a slight abuse of set-theoretic notation we define the class $\bigotimes_{i=1}^n \KK_i$ of $\Theta$-structures
as follows:
\begin{align*}
  \textstyle\bigotimes_{i=1}^n \KK_i = \big\{ \calA : \calA &\text{ is a $\Theta$-structure such that }\\
                                                  &\reduct{\calA}{\Sigma_i} \in \KK_i \text{ for all } i \in \{1, \ldots, n\}\big\}.
\end{align*}

\begin{THM}\cite[Corollary 2]{sokic-boron}\label{mp.thm.sokic}
  Let $\Sigma_1$, $\Sigma_2$, \ldots, $\Sigma_n$, $n \ge 2$, be pairwise disjoint sets of relational symbols,
  and for each $i \in \{1, 2, \ldots, n\}$ let $\KK_i$ be a class of $\Sigma_i$-structures having
  (HP), (JEP), (SAP) and the Ramsey property. Then $\bigotimes_{i=1}^n \KK_i$ has the Ramsey property.
\end{THM}

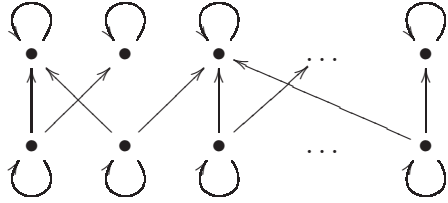
\begin{figure}[h]
$$
\xymatrix{
  & \\
  \bullet \ar@(ur,ul) & \bullet \ar@(ur,ul) & \bullet \ar@(ur,ul) & \ldots & \bullet \ar@(ur,ul) \\
  \bullet \ar@(dr,dl) \ar[u] \ar[ur] & \bullet \ar@(dr,dl) \ar[ur] \ar[ul] & \bullet \ar@(dr,dl) \ar[u] \ar[ur] & \ldots & \bullet \ar@(dr,dl) \ar[u] \ar[ull] \\
  & 
}
$$
\caption{A binary category}
\label{mp.fig.bc}
\end{figure}

In~\cite{masul-drp-perm} we devised a technique to transfer the Ramsey property from a category to
its (not necessarily full) subcategory, as follows.
Consider a finite, acyclic, bipartite digraph with loops 
where all the arrows go from one class of vertices into the other
and the out-degree of all the vertices in the first class is~2 (modulo loops), see Fig.~\ref{mp.fig.bc}.
Such a digraph can be thought of as a category (where the loops represent the identity morphisms), and will be referred to as
a \emph{binary category}.

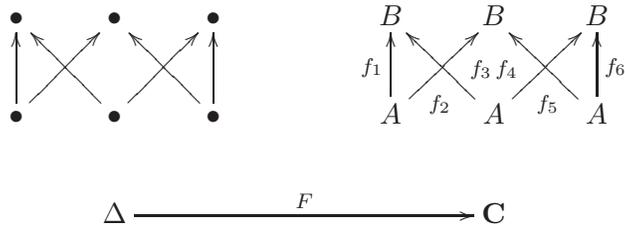
\begin{figure}[h]
  $$
  \xymatrix{
    \bullet & \bullet & \bullet
    & & B & B & B
  \\
    \bullet \ar[u] \ar[ur] & \bullet \ar[ur] \ar[ul] & \bullet \ar[ul] \ar[u]
    & & A \ar[u]^{f_1} \ar[ur]_(0.3){f_2} & A \ar[ur]^(0.3){f_4} \ar[ul]_(0.3){f_3} & A \ar[ul]^(0.3){f_5} \ar[u]_{f_6}
  \\
    & \Delta \ar[rrrr]^F  & & & & \CC  
  }
  $$
  \caption{A binary diagram in $\CC$ (of shape $\Delta$)}
  \label{catsap.fig.2}
\end{figure}

A \emph{binary diagram} in a category $\CC$ is a functor $F : \Delta \to \CC$ where $\Delta$ is a binary category,
$F$ takes the top row of $\Delta$ onto the same object, and takes the bottom row of $\Delta$ onto the same object,
Fig.~\ref{catsap.fig.2}.
A subcategory $\DD$ of a category $\CC$ is \emph{closed for binary diagrams} if every binary diagram
$F : \Delta \to \DD$ which is consistent in $\CC$ is also consistent in~$\DD$.

\begin{THM}\cite{masul-drp-perm}
  Let $\CC$ be a category such that $\hom_\CC(A, B)$ is finite for all $A, B \in \Ob(\CC)$
  and such that every morphism in $\CC$ is monic.
  Let $\DD$ be a (not necessarily full) subcategory of~$\CC$. If $\CC$ has the Ramsey property and $\DD$ is closed for binary diagrams,
  then $\DD$ has the Ramsey property.
\end{THM}

We will use the following special case of the previous theorem.

\begin{COR}\label{mp.cor.1}
  Let $\CC$ be a class of finite structures and $\KK$ a subclass of~$\CC$.
  If $\CC$ has the Ramsey property and $\KK$ is closed for binary diagrams (of structures and embeddings),
  then $\KK$ has the Ramsey property.
\end{COR}

\section{The proof}
\label{mp.sec.proof}

We are now ready to prove Theorem~\ref{mp.thm.MAIN}.

\bigskip

\begin{proof}
  Let  $\calT = (\{1, 2, \ldots, t\}, \Boxed\preccurlyeq)$, $t \in \NN$, be a template.
  Let $k_1$, $k_2$, \ldots, $k_m$ be all the isolated points in $\calT$, let $T' = T \setminus \{k_1$, $k_2$, \ldots, $k_m\}$ and let
  $\calT' = \restr{\calT}{T'}$. Clearly, here are no isolated points in $\calT'$, so every maximal element of $\calT'$ is above an element of $\calT'$.
  Let us fix a list
  $$
    (i_1, j_1), (i_2, j_2), \ldots, (i_s, j_s)
  $$
  of all the pairs $(i_\alpha, j_\alpha)$ of elements of $\calT'$ such that $i_\alpha \prec j_\alpha$
  in $\calT'$ and $j_\alpha$ is a maximal element in $\calT'$.
  
  Let $\overline\KK(\calT)$ be the class of structures defined as follows.
  For each
  $$
    (A, \Boxed{\le_1}, \ldots, \Boxed{\le_t}) \in \KK(\calT)
  $$
  the class $\overline\KK(\calT)$ contains
  $$
    (A, \Boxed{\le_{i_1}}, \Boxed{\le_{i_2}}, \ldots, \Boxed{\le_{i_s}},
    \Boxed{\le_{j_1}}, \Boxed{\le_{j_2}}, \ldots, \Boxed{\le_{j_s}},
    \Boxed{\le_{k_1}}, \Boxed{\le_{k_2}}, \ldots, \Boxed{\le_{k_m}})
  $$
  and these are the only structures in $\overline\KK(\calT)$.
  Clearly, $\overline\KK(\calT)$ and $\KK(\calT)$ are isomorphic as categories (where morphisms are embeddings),
  so it suffices to show that $\overline\KK(\calT)$ has the Ramsey property.

  It is easy to see that a structure $(A, \Boxed{\le_1}, \Boxed{\le_2}, \ldots, \Boxed{\le_{2s + m}})$
  with $2s + m$ binary relations on $A$ belongs to $\overline\KK(\calT)$ if and only if
  {\def\descriptionlabel#1{\hspace{\labelsep}#1}
  \begin{description}
  \item[(MP1)]
    $\le_\alpha$ are partial and $\le_{s + \beta}$ linear orders on $A$ for all $\alpha \in \{1, \ldots, s\}$
    and $\beta \in \{1, \ldots, s + m\}$;
  \item[(MP2)]
    $(\le_{\alpha}) \subseteq (\le_\beta)$ whenever
    $i_{\alpha} \preccurlyeq i_{\beta}$, for all $\alpha, \beta \in \{1, \ldots, s\}$;
  \item[(MP3)]
    $(\le_\alpha) \subseteq (\le_{s+\beta})$ whenever $i_{\alpha} \preccurlyeq j_{\beta}$, for all $\alpha, \beta \in \{1, \ldots, s\}$;
  \item[(MP4)]
    if $j_\alpha = j_\beta$ then $(\le_{s+\alpha}) = (\le_{s+\beta})$, for all $\alpha, \beta \in \{1, \ldots, s\}$.
  \end{description}}

  Let $\CC(s, m)$ be the class of structures of the form
  $(C, \Boxed{\le_1}, \Boxed{\le_2}, \ldots, \Boxed{\le_{2s+m}})$
  where $\le_\alpha$ are partial and $\le_{s + \beta}$ linear orders on $C$ for all $\alpha \in \{1, \ldots, s\}$ and
  $\beta \in \{1, \ldots, s + m\}$, and $\le_{s + \alpha}$ is a linear extension of $\le_\alpha$ for all~$\alpha \in \{1, 2, \ldots, s\}$.

  Let us show that $\CC(s, m)$ has the Ramsey property.
  Let $\Sigma_\alpha = \{\Boxed{\le_\alpha}, \Boxed{\le_{s + \alpha}}\}$ for $\alpha \in \{1, \ldots, s\}$ and
  $\Sigma_{s + \beta} = \{\Boxed{\le_{2s + \beta}}\}$ for $\beta \in \{1, \ldots, m\}$. Now, for $\alpha \in \{1, \ldots, s\}$ let
  $\KK_\alpha$ be the class $\EPosemb$ but over the signature $\Sigma_\alpha$, and for $\beta \in \{1, \ldots, m\}$ let
  $\KK_{s + \beta}$ be the class $\CHemb$ but over the signature $\Sigma_{s + \beta}$.
  It is easy to see that $\CC(s, m) = \bigotimes_{\alpha=1}^{s + m} \KK_\alpha$, whence follows that $\CC(s, m)$ has the
  Ramsey property by Theorem~\ref{mp.thm.sokic} and Examples~\ref{mp.ex.class-HP-JEP-SAP} and~\ref{mp.ex.class-RP}.

  As we have seen, $\overline\KK(\calT)$ is a subclass of a Ramsey class $\CC(s, m)$.
  By Corollary~\ref{mp.cor.1}, in order to show that $\overline\KK(\calT)$ has the Ramsey
  property it suffices to show that $\overline\KK(\calT)$ is closed for binary diagrams (of structures and embeddings)
  in~$\CC(s, m)$.

  Take any $\calA, \calB \in \overline\KK(\calT)$ and let $F : \Delta \to \overline\KK(\calT)$ be a binary diagram
  that takes the top row of $\Delta$ onto $\calB$ and the bottom row of $\Delta$ onto $\calA$.
  Assume that $F$ is consistent in $\CC(s, m)$ and
  let $\calC= (C, \Boxed{\le_1^\calC}, \Boxed{\le_2^\calC}, \ldots, \Boxed{\le_{2s+m}^\calC})$
  together with the embeddings
  $e_1, e_2, \ldots, e_n : \calB \hookrightarrow \calC$
  be a compatible cone in $\CC(s, m)$ over $F$.
  Define $\calD= (D, \Boxed{\le_1^\calD}, \Boxed{\le_2^\calD}, \ldots, \Boxed{\le_{2s+m}^\calD})$
  as follows. Let $D = e_1(B) \cup e_2(B) \cup \ldots \cup e_n(B)$.
  For every partial order $\sqsubseteq$ on $D$ there are many ways to choose a linear extension of $\sqsubseteq$.
  Let $\sqsubseteq_{lin}$ denote an arbitrary but fixed linear extension of $\sqsubseteq$ on $D$.
  Now, for each $\alpha \in \{1, \ldots, s\}$ and $\beta \in \{1, \ldots, s + m\}$ let
  $$
    \Boxed{\le_\alpha^\calD} =
    (
      \restr{\Boxed{\le_\alpha^\calC}}{e_1(B)} \cup
      \restr{\Boxed{\le_\alpha^\calC}}{e_2(B)} \cup \ldots \cup
      \restr{\Boxed{\le_\alpha^\calC}}{e_n(B)}
    )^+
  $$
  and
  $$
    \Boxed{\le_{s + \beta}^\calD} =
    ((
      \restr{\Boxed{\le_{s+\beta}^\calC}}{e_1(B)} \cup
      \restr{\Boxed{\le_{s+\beta}^\calC}}{e_2(B)} \cup \ldots \cup
      \restr{\Boxed{\le_{s+\beta}^\calC}}{e_n(B)}
    )^+)_{lin},
  $$
  where $^+$ denotes the transitive closure of a binary relation.
  Let us show that $\calD \in \overline\KK(\calT)$, or equivalently, that $\calD$ satisfies (MP1--4).
  {\def\descriptionlabel#1{\hspace{\labelsep}#1}
  \begin{description}
  \item[(MP1)]
    is obvious and follows directly from the definition of $\le_\alpha^\calD$ and $\le_{s+\beta}^\calD$.
  \item[(MP2)]
    is also easy to confirm. Assume that $i_{\alpha} \preccurlyeq i_{\beta}$ for some $\alpha, \beta \in \{1, 2, \ldots, s\}$.
    Then $(\Boxed{\le_\alpha^\calB}) \subseteq (\Boxed{\le_\beta^\calB})$ because $\calB \in \overline\KK(\calT)$.
    Since all the $e_i$'s are embeddings, $(\restr{\Boxed{\le_\alpha^\calC}}{e_1(B)}) \subseteq (\restr{\Boxed{\le_\beta^\calC}}{e_1(B)})$
    whence $(\Boxed{\le_\alpha^\calD}) \subseteq (\Boxed{\le_\beta^\calD})$.
  \item[(MP3) and (MP4)]
    follow by analogous arguments, having in mind that $\sqsubseteq_{lin}$ is an arbitrary but \emph{fixed} linear extension of~$\sqsubseteq$.
  \end{description}}

  Finally, define $f_1, f_2, \ldots, f_n : \calB \to \calD$ by $f_i(b) = e_i(b)$ for each $b \in B$ and $i \in \{1, \ldots, n\}$ and let us
  show that $f_i$'s are embeddings $\calB \hookrightarrow \calD$. Fix an $i \in \{1, \ldots, n\}$.
  Clearly, it suffices to prove that by taking transitive closures to obtain the partial and linear orders $\le_{\alpha}^\calD$
  and $\leqslant_{s + \beta}^\calD$ we do not change the orders restricted to $e_i(B)$.
  This is obvious for the linear orders, since the restriction to $e_i(B)$ of $\Boxed{\le_{s + \beta}^\calC}$ was already a linear order.
  As for the partial orders, notice that $(\Boxed{\le_\alpha^\calD}) \subseteq (\Boxed{\le_\alpha^\calC})$ for each $\alpha$,
  since $\Boxed{\le_\alpha^\calD}$ is a transitive closure of a subrelation of $\Boxed{\le_\alpha^\calC}$ and
  hence the restriction to $e_i(B)$ remains the same.

  Therefore, $\calD$ together with the embeddings $f_1, f_2, \ldots, f_n : \calB \hookrightarrow \calD$
  is a compatible cone in $\overline\KK(\calT)$ over $F$, which completes the proof.
\end{proof}

\section{Acknowledgements}

The authors would like to thank the two anonymous referees who shared their insights into the
historical development of identifying variuos Ramsey classes based on partial orders
and (often intertwined) paths to their solutions.

The second author gratefully acknowledges the support of the Grant No.\ 174019 of the Ministry of Education, Science and Technological Development
of the Republic of Serbia.


\begin{thebibliography}{99}\frenchspacing
\bibitem{AH}
  F.\ G.\ Abramson, L.\ A.\ Harrington.
  Models without indiscernibles.
  J.~Symbolic Logic 43 (1978), 572--600.

\bibitem{arman-rodl}
  A. Arman, V. R\"odl.
  Note on Ramsey theorem for posets with linear extensions.
  The Electronic Journal of Combinatorics 24.4 (2017), 4--36.

\bibitem{bodirsky1}
  M. Bodirsky.
  New Ramsey classes from old.
  The Electronic Journal of Combinatorics 21.2 (2014), 2--22.

\bibitem{bodirsky2}
  M. Bodirsky.
  Ramsey classes: Examples and constructions.
  Surveys in Combinatorics 2015, 424:1, 2015.

\bibitem{bottcher-foniok}
  J. B\"ottcher, J. Foniok.
  Ramsey Properties of Permutations.
  The Electronic Journal of Combinatorics 20.1 (2013) \#P2.

\bibitem{braunfeld}
  S. Braunfeld.
  Ramsey expansion of $\Lambda$-ultrametric spaces.
  Preprint, arXiv:1710.01193.

\bibitem{braunfeld-simon}
  S. Braunfeld, P. Simon.
  The classification of homogeneous finite-dimensional permutation structures.
  Preprint, arXiv:1807.07110.

\bibitem{cameron-perm}
  P. J. Cameron.
  Homogeneous permutations.
  The Electronic Journal of Combinatorics  9.2 (2002-3) \#R2.

\bibitem{fouche}
  W.\ L.\ Fouch\'e.
  Symmetry and the Ramsey degree of posets.
  15th British Combinatorial Conference (Stirling, 1995). Discrete Math.\ 167/168 (1997), 309--315.

\bibitem{GRS}
  R.\ L.\ Graham, B.\ L.\ Rothschild, J.\ H.\ Spencer.
  Ramsey Theory (2nd Ed).
  John Wiley \& Sons, 1990.

\bibitem{leeb-cat}
    K.\ Leeb.
    The categories of combinatorics.
    Combinatorial structures and their applications. Gordon and Breach, New York (1970).

\bibitem{masulovic-ramsey}
  D.\ Ma\v sulovi\'c, L.\ Scow.
  Categorical equivalence and the Ramsey property for finite powers of a primal algebra.
  Algebra Universalis 78 (2017), 159--179.

\bibitem{masul-preadj}
  D.\ Ma\v sulovi\'c.
  Pre-adjunctions and the Ramsey property.
  European Journal of Combinatorics 70 (2018), 268--283.

\bibitem{masul-drp-perm}
  D.\ Ma\v sulovi\'c.
  A Dual Ramsey Theorem for Permutations.
  Electronic Journal of Combinatorics 24(3) (2017), \#P3.39

\bibitem{masul-varieries-of-lattices}
  D.\ Ma\v sulovi\'c.
  The Ramsey and the ordering property for classes of lattices and semilattices.
  Order (to appear)

\bibitem{N1995}
  J.\ Ne\v set\v ril.
  Ramsey theory. In: R.\ L.\ Graham, M.\ Gr\"otschel and L.\ Lov\'asz, eds, Handbook of Combinatorics, Vol.~2,
  1331--1403, MIT Press, Cambridge, MA, USA, 1995.

\bibitem{nesetril}
  J.\ Ne\v set\v ril.
  Ramsey classes and homogeneous structures.
  Combinatorics, Probability and Computing 14 (2005), 171--189

\bibitem{Nesetril-Rodl}
  J.\ Ne\v set\v ril, V.\ R\"odl.
  Partitions of finite relational and set systems.
  J.\ Combin.\ Theory Ser.\ A 22 (1978), 289--312.

\bibitem{Nesetril-Rodl-1984}
  J.\ Ne\v set\v ril, V.\ R\"odl.
  Combinatorial partitions of finite posets and lattices -- Ramsey lattices.
  Algebra Universalis 19 (1984), 106--119.

\bibitem{Nesetril-Rodl-2016}
  J.\ Ne\v set\v ril, V.\ R\"odl.
  Ramsey Partial Orders from Acyclic Graphs.
  Order 35 (2018), 293--300.

\bibitem{PTW}
  M.\ Paoli, W.\ T.\ Trotter Jr., J.\ W.\ Walker.
  Graphs and orders in Ramsey theory and in dimension theory.
  Graphs and order (Banff, Alta., 1984), 351--394, NATO Adv.\ Sci.\ Inst.\ Ser.\ C Math.\ Phys.\ Sci., 147,
  Reidel, Dordrecht, 1985.

\bibitem{Ramsey}
  F.\ P.\ Ramsey.
  On a problem of formal logic.
  Proc.\ London Math.\ Soc.\ 30 (1930), 264--286.

\bibitem{sokic-phd}
  M.\ Soki\'c.
  Ramsey property of posets and related structures.
  PhD dissertation,
  University of Toronto, Canada, 2010.

\bibitem{sokic-order}
  M. Soki\'c.
  Ramsey properties of finite posets.
  Order 29 (2012), 1--30

\bibitem{sokic-boron}
  M.\ Soki\'c.
  Directed graphs and boron trees.
  J.\ Combin.\ Theory Ser.\ A 132 (2015), 142--171.

\bibitem{solecki-abstract}
  S.\ Solecki.
  Abstract approach to finite Ramsey theory and a self-dual Ramsey theorem.
  Adv.\ Math.\ 248 (2013) 1156--1198.

\bibitem{solecki-zhao}
  S.\ Solecki, M.\ Zhao.
  A Ramsey theorem for partial orders with linear extensions.
  European J.\ Comb.\ 60 (2017), 21--30.
\end{thebibliography}
\end{document}